\documentclass[12 pt]{amsart}
\usepackage{amssymb}
\usepackage{amsmath}
\usepackage{amscd}
\usepackage{graphicx}
\usepackage{latexsym}

\calclayout\evensidemargin .25in

\newcommand{\D}{\mathcal{D}}

\newcommand{\R}{\mathbb{R}}

\def\R{\mathbb{R}}

\newtheorem{theorem}{Theorem}[section]
\newtheorem{lemma}[theorem]{Lemma}
\newtheorem{proposition}[theorem]{Proposition}
\newtheorem{corollary}[theorem]{Corollary}
\newtheorem{definition}[theorem]{Definition}
\newtheorem{remark}[theorem]{Remark}

\def\dfn#1{{\em #1}}

\title[On the support genus of a contact structure] {On the support genus of a contact structure}
\author{Mehmet Firat Arikan}
\address{Department of Mathematics, MSU, East Lansing MI 48824, USA}
\email{arikanme@msu.edu}
\thanks{The author was partially supported by NSF Grant DMS0244622}

\begin{document}
\begin{abstract}
The algorithm given by Akbulut and Ozbagci constructs an explicit
open book decomposition on a contact three-manifold described by a
contact surgery on a link in the three-sphere. In this article, we
will improve this algorithm by using Giroux's contact cell
decomposition process. In particular, our algorithm gives a better
upper bound for the recently defined ``minimal supporting genus
invariant'' of contact structures.
\end{abstract}

\maketitle

\section{Introduction}

\medskip \noindent
Let $(M,\xi)$ be a closed oriented contact 3-manifold, and let
$(\Sigma,h)$ be an open book (decomposition) of $M$ which is
compatible with the contact structure $\xi$ (sometimes we also say
that $(\Sigma,h)$ supports $\xi$). Based on the correspondence
theorem (see Theorem \ref{Giroux}) between contact structures and
their supporting open books, the topological invariant $sg(\xi)$
was defined in \cite{EO}. More precisely, we have
\[
sg(\xi)=\min \{ \, \,g(\Sigma) \, \, \vert \, \,(\Sigma,h) \text{ an
open book decomposition supporting } \xi\}
\]
called \dfn{supporting genus} of $\xi$. There are some partial
results for this invariant. For instance, we have:

\begin{theorem} [\cite{Et1}] \label{sg=0ifovertwisted}
If $(M,\xi)$ is overtwisted, then $sg(\xi)=0.$
\end{theorem}

\noindent Unlike the overtwisted case, there is not much known yet
for $sg(\xi)$ when $\xi$ is tight. On the other hand, if we,
furthermore, require that $\xi$ is Stein fillable, then an
algorithm to find an open book supporting $\xi$ was given in
\cite{AO}. Although their construction is explicit, the pages of
the resulting open books arise as Seifert surfaces of torus knots
or links, and so this algorithm is far from even approximating the
numbers $sg(\xi)$. In \cite{St}, the same algorithm was
generalized to the case where $\xi$ need not to be Stein fillable
(or even tight), but the pages are still of large genera.

\medskip \noindent This article is organized as follows: After the 
preliminaries (Section \ref{preliminaries}), in Section \ref{Algorithm} 
we will present an explicit construction 
of a supporting open book (with considerably less genus) for a
given contact surgery diagram of any contact structure $\xi$. Of
course, because of Theorem \ref{sg=0ifovertwisted}, our algorithm
makes more sense for the tight structures than the overtwisted
ones. Moreover, it depends on a choice of the contact surgery
diagram describing $\xi$. Nevertheless, it gives better and more
reasonable upper bound for $sg(\xi)$ (when $\xi$ is tight) as we
will see from our examples in Section \ref{examples}.

\medskip \noindent Let $L$ be any Legendrian link given in
$(\mathbb{R}^3,\xi_0=ker(\alpha_0= dz+xdy))\subset(S^3,\xi_{st})$.
$L$ can be represented by a special diagram $\D$ called \emph{a
square bridge diagram} of $L$ (see \cite{Ly}). We will consider
$\D$ as an abstract diagram such that
\begin{enumerate}
\item $\D$ consists of horizontal line segments $h_1, ..., h_p$, and vertical
line segments $v_1,..., v_q$ for some integers $p\geq2$, $q\geq2$,
\item there is no collinearity in $\{h_1, \dots, h_p \}$, and in $\{v_1, \dots,
v_q \}$.
\item each $h_i$ (resp., each $v_j$) intersects two vertical
(resp., horizontal) line segments of $\D$ at its two endpoints
(called \emph{corners} of $\D$), and
\item any interior intersection (called \emph{junction} of $\D$)
 is understood to be a virtual crossing of $\D$
where the horizontal line segment is passing over the vertical one.
\end{enumerate}

\begin{figure}[ht]
  \begin{center}
   \includegraphics{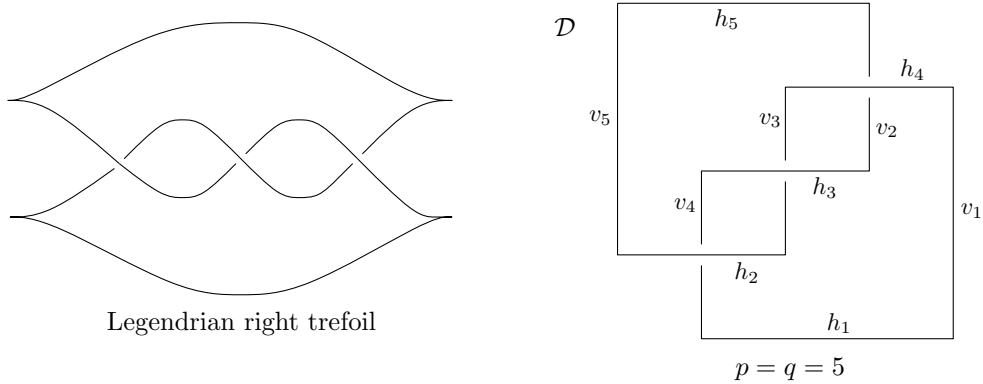}
   \caption{The square bridge diagram $\D$ for the Legendrian right trefoil}
  \label{trefoil_0}
    \end{center}
\end{figure}

\noindent We depict Legendrian right trefoil and the corresponding $\D$ in Figure
\ref{trefoil_0}.

\noindent Clearly, for any front projection of a Legendrian link,
we can associate a square bridge diagram $\D$. Using such a
diagram $\D$, the following two facts were first proved in
\cite{AO}, and later made more explicit in \cite{Pl}. Below
versions are from the latter:

\begin{lemma}\label{page} Given a Legendrian link $L$ in $(\R^3, \xi_0)$,
there exists a torus link $T_{p,q}$ (with $p$ and $q$ as above)
transverse to $\xi_0$ such that its Seifert surface $F_{p,q}$
contains $L$, $d\alpha_0$ is an area form on $F_{p,q}$, and $L$
does not separate $F_{p,q}$.
\end{lemma}

\begin{proposition}\label{linkinbook} Given $L$ and $F_{p,q}$ as above,
there exist an open book decomposition of $S^3$ with page $F_{p,q}$
such that:
\begin{enumerate}
\item
 the induced contact structure $\xi$ is isotopic
to $\xi_0$;
\item
   the link $L$ is contained in one of the page $F_{p,q}$, and does not separate it;
\item
$L$ is Legendrian with respect to $\xi$;
\item
there exist an isotopy which fixes $L$ and takes $\xi$ to $\xi_0$,
so the Legendrian type of the link is the same with respect to $\xi$
and $\xi_0$;
\item the framing of $L$ given by the page $F_{p,q}$ of the open book is the
same as the contact framing.
\end{enumerate}
\end{proposition}

\noindent Being a Seifert surface of a torus link, $F_{p,q}$ is of
large genera. In Section \ref{Algorithm}, we will construct
another open book $\mathcal{OB}$ supporting $(S^3,\xi_{st})$ such
that its page $F$ arises as a subsurface of $F_{p,q}$ (with
considerably less genera), and given Legendrian link $L$ sits on
$F$ as how it sits on the page $F_{p,q}$ of the construction used
in \cite{AO} and \cite{Pl}. The page $F$ of the open book
$\mathcal{OB}$ will arise as the ribbon of the 1-skeleton of an
appropriate contact cell decomposition for $(S^3,\xi_{st})$. As in
\cite{Pl}, our construction will keep the given link $L$
Legendrian with respect to the standard contact structure
$\xi_{st}$. Our main theorem is:

\begin{theorem}\label{existence_of_CCD} Given $L$ and $F_{p,q}$ as above,
there exists a contact cell decomposition $\Delta$ of
$(S^3,\xi_{st})$ such that
\begin{enumerate}
\item $L$ is contained in the Legendrian 1-skeleton $G$ of $\Delta$,
\item The ribbon $F$ of the 1-skeleton $G$ is a subsurface of
$F_{p,q}$ ($p$ and $q$ as above),
\item The framing of $L$ coming from $F$ is equal to its contact framing
$tb(L)$, and
\item If $p>3$ and $q>3$, then the
genus $g(F)$ of $F$ is strictly less than the genus $g(F_{p,q})$ of
$F_{p,q}$.
\end{enumerate}
\end{theorem}

\noindent As an immediate consequence (see Corollary
\ref{main_corollary}), we get an explicit description of an open
book supporting $(S^3,\xi)$ whose page $F$ contains $L$ with the
correct framing. Therefore, if $(M^{\pm},\xi^{\pm})$ is given by
contact ($\pm 1$)-surgery on $L$ (such a surgery diagram exists
for any closed contact 3-manifold by Theorem \ref{DingGeigesmain}),
we get an open book supporting $\xi^{\pm}$ with page $F$ by
Theorem \ref{surgeryonbook}. Hence, $g(F)$ improves the upper bound
for $sg(\xi)$ as $g(F) < g(F_{p,q})$ (for $p>3, \; q>3$).
It will be clear from our examples in Section \ref{examples}
that this is indeed a good improvement.

\noindent {\em Acknowledgments.\/} The author would like to thank
Selman Akbulut, Selahi Durusoy, Cagri Karakurt, and Burak Ozbagci
for their helpful conversations and comments  on the draft of this
paper.


\section{Preliminaries} \label{preliminaries}
\subsection{Contact structures and Open book decompositions}

\noindent A $1$-form $\alpha \in \Omega^1(M)$ on a $3$-dimensional
oriented manifold $M$ is called a \dfn{contact form} if it satisfies
$\alpha \wedge d\alpha \neq 0$. An \dfn{oriented contact structure}
on $M$ is then a hyperplane field $\xi$ which can be globally
written as the kernel of a contact $1$-form $\alpha$. We will always
assume that $\xi $ is a \dfn{positive} contact structure, that is,
$\alpha \wedge d\alpha > 0$. Note that this is equivalent to asking
that $d\alpha$ be positive definite on the plane field $\xi$, ie.,
$d\alpha |_{\xi}>0$. Two contact structures $\xi_0, \xi_1$ on a
$3$-manifold are said to be \emph{isotopic} if there exists a
1-parameter family $\xi_t$ ($0\leq t\leq 1$) of contact structures
joining them. We say that two contact $3$-manifolds $(M_1,\xi_1)$
and $(M_2,\xi_2)$ are \emph{contactomorphic} if there exists a
diffeomorphism $f:M_1\longrightarrow M_2$ such that
$f_\ast(\xi_1)=\xi_2$. Note that isotopic contact structures give
contactomorphic contact manifolds by Gray's Theorem. Any contact
$3$-manifold is locally contactomorphic to $(\mathbb{R}^3,\xi_0)$
where \emph{standard contact structure} $\xi_0$ on $\mathbb{R}^3$
with coordinates $(x,y,z)$ is given as the kernel of
$\alpha_0=dz+xdy$. The standard contact structure $\xi_{st}$ on the
$3$-sphere $S^3=\{(r_1,r_2,\theta_1,\theta_2) : r_1^2+r_2^2=1\}
\subset \mathbb{C}^2$ is given as the kernel of
$\alpha_{st}=r_1^2d\theta_1+r_2^2d\theta_2$. One basic fact is that
$(\mathbb{R}^3, \xi_0)$ is contactomorphic to $(S^3 \setminus \{pt
\},\xi_{st})$. For more details on contact geometry, we refer the
reader to \cite{Ge}, \cite{Et3}.

\medskip \noindent
\noindent An \dfn{open book decomposition} of a closed $3$-manifold
$M$ is a pair $(L, f)$ where $L$ is an oriented link in $M$, called
the \dfn{binding}, and $f: M\setminus L \to S^1$ is a fibration such
that $f^{-1}(t)$ is the interior of a compact oriented surface
$\Sigma_t \subset M$ and $\partial \Sigma_t=L$ for all $t \in S^1$.
The surface $\Sigma=\Sigma_t$, for any $t$, is called the \dfn{page}
of the open book. The \dfn{monodromy} of an open book $(L,f)$ is
given by the return map of a flow transverse to the pages (all
diffeomorphic to $\Sigma$) and meridional near the binding, which is
an element $h \in Aut(\Sigma,\partial \Sigma)$, the group of
(isotopy classes of) diffeomorphisms of $\Sigma$ which restrict to
the identity on $\partial \Sigma$ . The group $Aut(\Sigma,\partial
\Sigma)$ is also said to be the mapping class group of $\Sigma$, and
denoted by $\Gamma(\Sigma)$.

\medskip \noindent
An open book can also be described as follows. First consider the
mapping torus $$\Sigma(h)= [0,1]\times \Sigma/(1,x)\sim (0, h(x))$$
where $\Sigma$ is a compact oriented surface with $n=|\partial
\Sigma|$ boundary components and $h$ is an element of
$Aut(\Sigma,\partial \Sigma)$ as above. Since $h$ is the identity
map on $\partial \Sigma$, the boundary $\partial \Sigma(h)$ of the
mapping torus $\Sigma(h)$ can be canonically identified with $n$
copies of $T^2 = S^1 \times S^1$, where the first $S^1$ factor is
identified with $[0,1] / (0\sim 1)$ and the second one comes from a
component of $\partial \Sigma$. Now we glue in $n$ copies of
$D^2\times S^1$ to cap off $\Sigma(h)$ so that $\partial D^2$ is
identified with $S^1 = [0,1] / (0\sim 1)$ and the $S^1$ factor in
$D^2 \times S^1$ is identified with a boundary component of
$\partial \Sigma$. Thus we get a closed $3$-manifold
$$M=M_{(\Sigma,h)}:= \Sigma(h) \cup_{n} D^2 \times S^1 $$
\emph{equipped with} an open book decomposition $(\Sigma,h)$ whose
binding is the union of the core circles in the $D^2 \times S^1$'s
that we glue to $\Sigma(h)$ to obtain $M$.
To summarize, an element $h \in Aut(\Sigma,\partial \Sigma)$
determines a $3$-manifold $M=M_{(\Sigma,h)}$ together with an
``abstract" open book decomposition $(\Sigma,h)$ on it. For furher
details on these subjects, see \cite{Gd}, and \cite{Et2}.


\subsection{Legendrian Knots and Contact Surgery }
A {\em Legendrian knot} $K$ in a contact $3$-manifold $(M,\xi )$
is a knot that is everywhere tangent to $\xi$. Any Legendrian knot
comes with a canonical {\em contact framing} (or
\emph{Thurston-Bennequin framing}), which is defined by a vector
field along $K$ that is transverse to $\xi$. If $K$ is null-homologous,
then this framing can be given by an integer $tb(K)$, called  \emph{Thurston-Bennequin number}.
For any Legendrian knot $K$ in $(\mathbb{R}^3, \xi_0)$, the number $tb(K)$ can be computed as
$$tb(K)=bb(K) - \#  \mbox{left cusps of K}$$
where $bb(K)$ is the blackboard framing of $K$.

\medskip \noindent We call $(M,\xi )$ (or just $\xi$) {\em overtwisted} if it contains an embedded disc
$D \approx D^2\subset M$ with boundary $\partial D \approx S^1$ a
Legendrian knot whose contact framing equals the framing it
receives from the disc $D$. If no such disc exists, the contact
structure $\xi$ is called {\em tight}.

\medskip \noindent For any $p,q \in \mathbb{Z}$, a \emph{contact
$(r)$-surgery} ($r=p/q$) along a Legendrian knot $K$ in a contact
manifold $(M,\xi )$ was first described in \cite{DG1}. It is
defined to be a special kind of a topological surgery, where
surgery coefficient $r\in{\mathbb Q}\cup{\infty}$ measured
relative to the contact framing of $K$. For $r\neq 0$, a contact
structure on the surgeried manifold
\[ (M-\nu K)\cup (S^1\times D^2),\]
($\nu K$ denotes a tubular neighborhood of $K$) is defined by
requiring this contact structure to coincide with $\xi$ on $Y-\nu
K$ and its extension over $S^1\times D^2$ to be tight on (glued
in) solid torus $S^1\times D^2$. Such an extension uniquely exists
(up to isotopy) for $r=1/k$ with $k\in \mathbb{Z}$ (see
\cite{Ho}). In particular, a contact $(\pm1)$-surgery along a
Legendrian knot $K$ on a contact manifold $(M,\xi)$ determines a
unique (up to contactomorphism) surgered contact manifold which
will be denoted by $(M,\xi)_{(K,\pm1)}$.

\medskip \noindent The most general result along these lines is:
\begin{theorem} [\cite{DG1}] \label{DingGeigesmain}
Every (closed, orientable) contact $3$-manifold $(M, \xi )$ can be
obtained via contact $(\pm 1)$-surgery on a Legendrian link in
$(S^3, \xi _{st})$.
\end{theorem}

\medskip \noindent
Any closed contact $3$-manifold $(M,\xi)$ can be described by a
\emph{contact surgery diagram}. Such a diagram consists of a front
projection (onto the $yz$-plane) of a Legendrian link drawn in
$(\mathbb{R}^3, \xi_0) \subset (S^3,\xi_{st})$ with contact
surgery coefficient on each link component. Theorem
\ref{DingGeigesmain} implies that there is a contact surgery
diagram for $(M,\xi)$ such that the contact surgery coefficient of
any Legendrian knot in the diagram is $\pm1$. For more details see
\cite{Gm} and \cite{OS}.


\subsection{Compatibility and Stabilization} A contact structure
$\xi$ on a $3$-manifold $M$ is said to be \dfn{supported by an open
book} $(L,f)$ if $\xi$ is isotopic to a contact structure given by a
$1$-form $\alpha$ such that
\begin{enumerate}
\item $d\alpha$ is a positive area form on each page $\Sigma \approx f^{-1}($pt$)$ of the open book and
\item $\alpha>0$ on $L$ (Recall that $L$ and the pages are oriented.)
\end{enumerate}

\noindent When this holds, we also say that the open book $(L,f)$
is \dfn{compatible with the contact structure} $\xi$ on $M$.
Geometrically, compatibility means that $\xi$ can be isotoped to
be arbitrarily close (as oriented plane fields), on compact
subsets of the pages, to the tangent planes to the pages of the
open book in such a way that after some point in the isotopy the
contact planes are transverse to $L$ and transverse to the pages
of the open book in a fixed neighborhood of $L$.

\begin{definition} \label{stabilization}
A positive (resp., negative) stabilization $S^{+}_{K}(\Sigma,h)$
(resp., $S^{-}_K(\Sigma,h)$) of an abstract open book $(\Sigma,h)$
is the open book
\begin{enumerate}
\item with page $\Sigma'=\Sigma \cup \text{ 1-handle}$ and
\item monodromy $h'=h \circ D_K$ (resp., $h'=h \circ D_K^{-1}$)
where $D_K$ is a right-handed Dehn twist along a curve $K$ in
$\Sigma'$ that intersects the co-core of the 1-handle exactly once.
\end{enumerate}
\end{definition}

\medskip \noindent
Based on the result of Thurston and Winkelnkemper \cite{TW},
Giroux proved the following theorem which strengthened the link
between open books and contact structures.

\begin{theorem} [\cite{Gi}] \label{Giroux}
Let $M$ be a closed oriented $3$-manifold. Then there is a
one-to-one correspondence between oriented contact structures on $M$
up to isotopy and open book decompositions of $M$ up to positive
stabilizations: Two contact structures supported by the same open
book are isotopic, and two open books supporting the same contact
structure have a common positive stabilization.
\end{theorem}

\noindent For a given fixed open book $(\Sigma,h)$ of a $3$-manifold
$M$, there exists a unique compatible contact structure up to
isotopy on $M=M_{(\Sigma,h)}$ by Theorem \ref{Giroux}. We will
denote this contact structure by $\xi_{(\Sigma,h)}$. Therefore, an
open book $(\Sigma,h)$ determines a unique contact manifold
$(M_{(\Sigma,h)},\xi_{(\Sigma,h)})$ up to contactomorphism.

\medskip \noindent Taking a positive stabilization of an open book
$(\Sigma,h)$ is actually taking a special Murasugi sum of
$(\Sigma,h)$ with $(H^{+},D_c)$ where $H^+$ is the positive Hopf
band, and $c$ is the core circle in $H^{+}$. Taking a Murasugi sum
of two open books corresponds to taking the connect sum of
$3$-manifolds associated to the open books. For the precise
statements of these facts, and a proof of the following theorem,
we refer the reader to \cite{Gd}, \cite{Et2}.
\begin{theorem} \label{connectsummingwithS^3}
$(M_{S^{+}_{K}(\Sigma,h)},\xi_{S^{+}_{K}(\Sigma,h)}) \cong
(M_{(\Sigma,h)},\xi_{(\Sigma,h)}) \# (S^3,\xi_{st}) \cong
(M_{(\Sigma,h)},\xi_{(S,h)}).$
\end{theorem}


\subsection{Monodromy and Surgery Diagrams}
Given a contact surgery diagram for a closed contact $3$-manifold
$(M,\xi)$, we want to construct an open book compatible with
$\xi$. One implication of Theorem \ref{DingGeigesmain} is that one
can obtain such a compatible open book by starting with a
compatible open book of $(S^3,\xi_{st})$, and then interpreting
the effects of surgeries (yielding $(M,\xi)\,$) in terms of open
books. However, we first have to realize each surgery curve (in
the given surgery diagram of $(M,\xi)\,$) as a Legendrian curve
sitting on a page of some open book supporting $(S^3,\xi_{st})$.
We refer the reader to Section 5 in \cite{Et2} for a proof of the
following theorem.

\begin{theorem} \label{surgeryonbook}
Let $(\Sigma,h)$ be an open book supporting the contact manifold
$(M,\xi).$ If $K$ is a Legendrian knot on the page $\Sigma$ of the
open book, then
\[
(M,\xi)_{(K,\; \pm 1)}= (M_{(\Sigma,\; h\circ D_K^{\mp})},
\xi_{(\Sigma,\; h \circ D_K^{\mp})}).
\]
\end{theorem}


\subsection{Contact Cell Decompositions and Convex Surfaces}
The exploration of contact cell decompositions in the study of
open books was originally initiated by Gabai \cite{Ga}, and then
developed by Giroux \cite{Gi}. We want to give several definitions and facts carefully.

\medskip \noindent Let $(M,\xi)$ be any contact 3-manifold,
and $K \subset M$ be a Legendrian knot. The \emph{twisting number}
$tw(K,Fr)$ of $K$ with respect to a given framing $Fr$ is defined
to be the number of counterclockwise $2\pi$ twists of $\xi$ along
$K$, relative to $Fr$. In particular, if $K$ sits on a surface
$\Sigma \subset M$, and $Fr_{\Sigma}$ is the surface framing of
$K$ given by $\Sigma$, then we write $tw(K, \Sigma)$ for $tw(K,
Fr_{\Sigma})$. If  $K=\partial\Sigma$, then we have $tw(K,
\Sigma)=tb(K)$ (by the definition of $\, tb$).

\begin{definition} \label{Cont_Cell_Dec}
A contact cell decomposition of a contact $3-$manifold $(M,\xi)$
is a finite CW-decomposition of $M$ such that
\begin{enumerate}
\item[(1)] the 1-skeleton is a Legendrian graph,
\item[(2)] each 2-cell $D$ satisfies $tw(\partial D, D)=-1,$ and
\item[(3)] $\xi$ is tight when restricted to each 3-cell.
\end{enumerate}
\end{definition}

\begin{definition}
Given any Legendrian graph $G$ in $(M,\xi)$, the ribbon of $G$ is
a compact surface $R=R_G$ satisfying
\begin{enumerate}
\item $R$ retracts onto $G,$
\item $T_p R=\xi_p$ for all $p\in G,$
\item $T_p R\not = \xi_p$ for all $p\in R\setminus G.$
\end{enumerate}
\end{definition}

\noindent For a proof of the following lemma we refer the reader to
\cite{Gd} and \cite{Et2}.

\begin{lemma} \label{RibbonIsPage}
Given a closed contact $3-$manifold $(M,\xi)$,
the ribbon of the 1-skeleton of any contact cell decomposition is a
page of an open book supporting $\xi$.
\end{lemma}

\noindent The following lemma will be used in the next section.

\begin{lemma} \label{ExtendingDecomposition}
Let $\Delta$ be a contact cell decomposition of a closed contact
3-manifold $(M,\xi)$ with the $1-$skeleton $G$. Let $U$ be a
3-cell in $\Delta$. Consider two Legendrian arcs $I \subset \partial U$ and $J \subset U$ such that
\begin{enumerate}
\item $I \subset G$,
\item $J \cap \partial U=\partial J=\partial I$,
\item $C=I \cup_{\partial} J$ is a Legendrian unknot with $tb(C)=-1$.
\end{enumerate}
Set $G'=G \cup J$. Then there exists another contact cell decomposition
$\Delta'$ of $(M,\xi)$ such that $G'$ is the 1-skeleton of
$\Delta'$
\end{lemma}

\begin{proof} The interior of the $3-$cell $U$ is contactomorphic to
$(\mathbb{R}^3,\xi_0)$. Therefore, there exists an embedded disk
$D$ in $U$ such that $\partial D=C$ and $int(D) \subset int(U)$ as
depicted in Figure \ref{J_contained_in_U}(a). We have $tw(\partial
D,D)= -1$ since $tb(C)=-1$. As we are working in
$(\mathbb{R}^3,\xi_0)$, there exist two $C^{\infty}$-small
perturbations of $D$ fixing $\partial D=C$ such that perturbed
disks intersect each other only along their common boundary $C$.
In other words, we can find two isotopies $H_1, H_2 : [0,1]\times
D \longrightarrow U$ such that for each $i=1,2$ we have
\begin{enumerate}
\item $H_i(t,.)$ fixes $\partial D=C$ pointwise for all $t \in
[0,1]$,
\item $H_i(0,D)=Id_D$ where $Id_D$ is the identity map on $D$,
\item $H_i(1,D)=D_i$ where each $D_i$ is an embedded disk in $U$ with $int(D_i) \subset
int(U)$,
\item $D \cap D_1 \cap D_2=C$ (see Figure \ref{J_contained_in_U}(b)).
\end{enumerate}

\begin{figure}[ht]
  \begin{center}
   \includegraphics{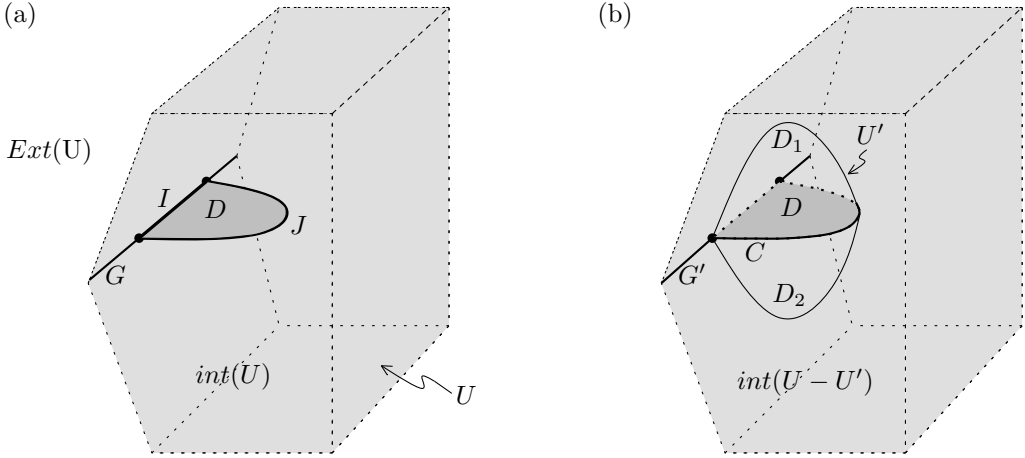}
   \caption{Constructing a new contact cell decomposition}
  \label{J_contained_in_U}
    \end{center}
\end{figure}

\noindent Note that $tw(\partial D_i,D_i)=tw(C,D_i)=-1$ for
$i=1,2$. This holds because each $D_i$ is a small perturbation of
$D$, so the number of counterclockwise twists of $\xi$ (along $K$)
relative to $Fr_{D_i}$ is equal to the one relative to $Fr_{D}$.

\medskip \noindent Next, we introduce $G'=G \cup J$ as the
1-skeleton of the new contact cell decomposition $\Delta'$. In
$M-int(U)$, we define the 2- and 3- skeletons of $\Delta'$ to be
those of $\Delta$ . However, we change the cell structure of
$int(U)$ as follows: We add 2-cells $D_1, D_2$ to the 2-skeleton
of $\Delta'$ (note that they both satisfy the twisting condition
in Definition \ref{Cont_Cell_Dec}). Consider the 2-sphere $S=D_1
\cup D_2$ where the union is taken along the common boundary $C$.
Let $U'$ be the 3-ball with $\partial U'=S$. Note that $\xi|_{U'}$
is tight as $U' \subset U$ and $\xi|U$ is tight. We add $U'$ and
$U-U'$ to the 3-skeleton of $\Delta'$ (note that $U-U'$ can be
considered as a 3-cell because observe that $int(U-U')$ is
homeomorphic to the interior of a 3-ball as in Figure
\ref{J_contained_in_U}(b)). Hence, we established another contact
cell decomposition of $(M,\xi)$ whose 1-skeleton is $G'=G\cup J$.
(Equivalently, by Theorem \ref{connectsummingwithS^3}, we are
taking the connect sum of $(M,\xi)$ with $(S^3,\xi_{st})$ along
$U'$.)
\end{proof}


\section{The Algorithm} \label{Algorithm}
\subsection{Proof of Theorem \ref{existence_of_CCD}}
\begin{proof} By translating $L$ in $(\R^3, \xi_0)$ if necessary
(without changing its contact type), we can assume that the front
projection of $L$ onto the $yz$-plane lying in the second quadrant
$\{\; (y,z) \;| \; y<0, \; z>0 \}$. After an appropriate
Legendrian isotopy, we can assume that $L$ consists of the line
segments contained in the lines
\begin{center} $k_i=\{x=1, z=-y+a_i\}$, $i=1,\dots, p$, \\
$l_j=\{x=-1, z=y+b_j\}$, $j=1,\dots, q$
\end{center}
for some $a_1<a_2< \dots <a_p$, $0<b_1<b_2< \dots <b_q$, and also
the line segments (parallel to the $x$-axis) joining certain
$k_i$'s to certain $l_j$'s. In this representation, $L$ seems to
have corners. However, any corner of $L$ can be made smooth by a
Legendrian isotopy changing only a very small neighborhood of that
corner.

\medskip \noindent Let $\pi:\mathbb{R}^3\longrightarrow
\mathbb{R}^2$ be the projection onto the $yz$-plane. Then we
obtain the square bridge diagram $D=\pi(L)$ of $L$ such that $D$
consists of the line segments
\begin{eqnarray*}
h_i \subset \pi(k_i) &=& \{x=0, z=-y+a_i\}, \quad i=1,\dots, p, \\
v_j \subset \pi(l_j) &=& \{x=0, z=y+b_j\}, \quad j=1,\dots, q.
\end{eqnarray*}
Notice that $D$ bounds a polygonal region $P$ in the second
quadrant of the $yz$-plane, and divides it into finitely many
polygonal subregions $P_1, \dots, P_m$ ( see Figure
\ref{trefoil_1}-(a) ).

\medskip \noindent Throughout the proof, we will assume that the link $L$ is
not split (that is, the region $P$ has only one connected
component). Such a restriction on $L$ will not affect the generality
of our construction (see Remark \ref{increasing_efficiency}).

\begin{figure}[ht]
  \begin{center}
   \includegraphics{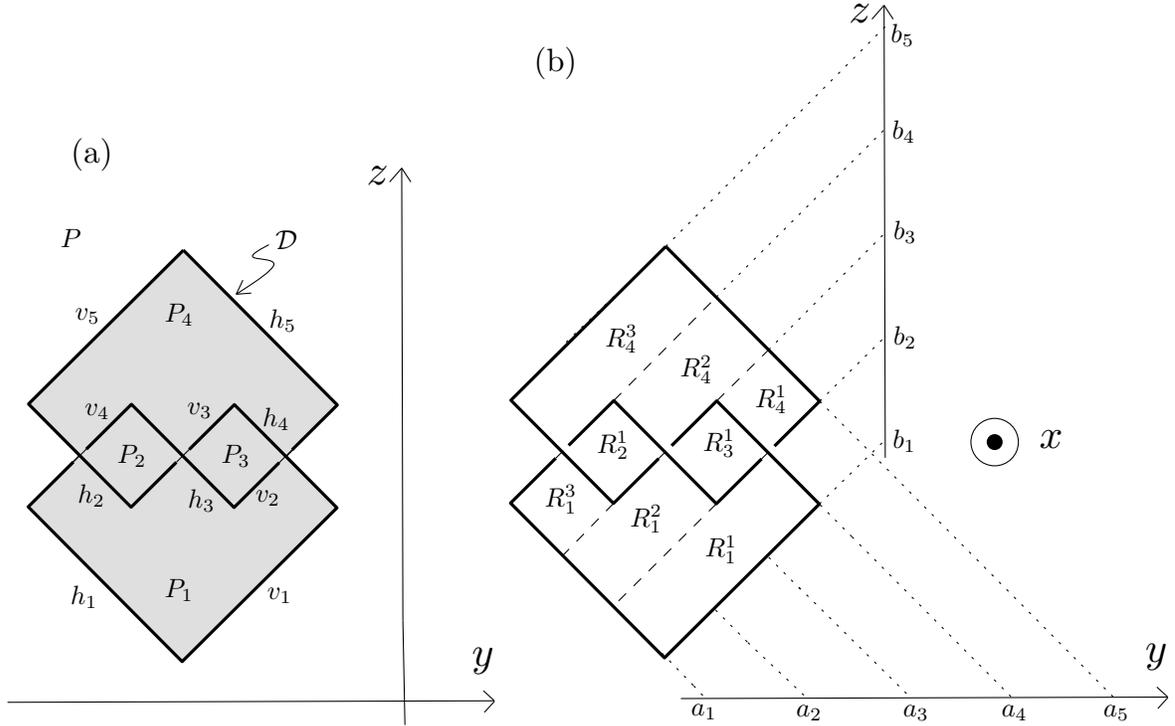}
   \caption{The region $P$ for right trefoil knot and its division into rectangles}
  \label{trefoil_1}
    \end{center}
\end{figure}

\noindent Now we decompose $P$ into finite number of ordered
rectangular subregions as follows: The collection $\{\pi(l_j)\;
|\; j=1,\dots, q \}$ cuts each $P_k$ into finitely many
rectangular regions $R_k^1, \dots, R_k^{m_k}$. Consider the set
$\mathfrak{P}$ of all such rectangles in $P$. That is, we define
\begin{center}
$\mathfrak{P} \doteq \{ \; R_k^{l} \;| \; k=1, \dots, m, \quad l=1,
\dots, m_k\}$.
\end{center}
Clearly $\mathfrak{P}$ decomposes $P$ into rectangular regions (
see Figure \ref{trefoil_1}-(b) ). The boundary of an arbitrary
element $R_k^{l}$ in $\mathfrak{P}$ consists of four edges: Two of
them are the subsets of the lines $\pi(l_{j(k,l)})$,
$\pi(l_{j(k,l)+1})$, and the other two are the subsets of the line
segments $h_{i_1(k,l)}$, $h_{i_2(k,l)}$ where $1\leq i_1(k,l) <
i_2(k,l) \leq p$ and $1\leq j(k,l) < j(k,l)+1 \leq q$ (see Figure
\ref{typical_rectangle}).

\begin{figure}[ht]
  \begin{center}
   \includegraphics{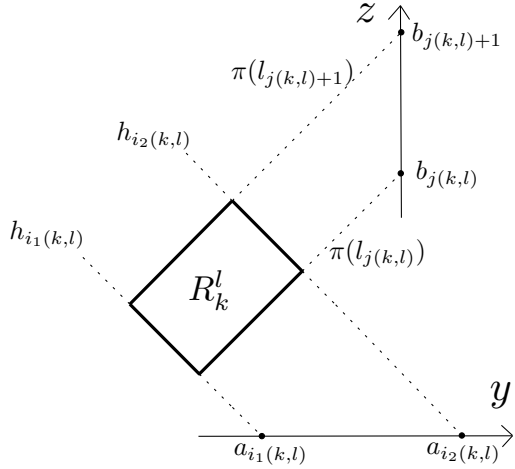}
   \caption{Arbitrary element $R_k^{l}$ in $\mathfrak{P}$}
  \label{typical_rectangle}
    \end{center}
\end{figure}

\noindent Since the region $P$ has one connected component, the
following holds for the set $\mathfrak{P}$:

\medskip \noindent
($\star$) Any element of $\mathfrak{P}$ has at least one common
vertex with some other element of $\mathfrak{P}$.
\medskip

\noindent By ($\star$), we can rename the elements of $\mathfrak{P}$
by putting some order on them so that any element of $\mathfrak{P}$
has at least one vertex in common with the union of all rectangles
coming before itself with respect to the chosen order. More
precisely, we can write
$$\mathfrak{P}=\{ \; R_k \; | \; k=1, \dots, N \}$$
($N$ is the total number of rectangles in $\mathfrak{P}$) such
that each $R_k$ has at least one vertex in common with the union
$R_1 \cup \dots \cup R_{k-1}$.

\medskip \noindent Equivalently, we can construct the polygonal
region $P$ by introducing the building rectangles ($R_k$'s) one by
one in the order given by the index set $\{ 1, 2, \dots, N\}$. In
particular, this eliminates one of the indexes, i.e., we can use
$R_k$'s instead of $R_k^l$'s. In Figure \ref{trefoil_2}, how we
build $P$ is depicted for the right trefoil knot (compare it with
the previous picture given for $P$ in Figure \ref{trefoil_1}-(b)).

\begin{figure}[ht]
  \begin{center}
   \includegraphics{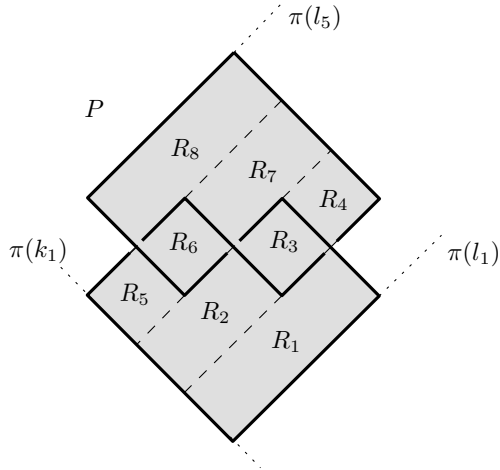}
   \caption{The region $P$ for right trefoil knot}
  \label{trefoil_2}
    \end{center}
\end{figure}

\medskip \noindent Using the representation $P=R_1\cup R_2 \cup
\dots \cup R_N$, we will construct the contact cell decomposition
(CCD) $\Delta$. Consider the following infinite strips which are
parallel to the $x$-axis (they can be considered as the unions of
``small'' contact planes along $k_i$'s and $l_j$'s):
\begin{eqnarray*}
S^+_i&=&\{ 1-\epsilon \leq x\leq 1+\epsilon, \; z=y+a_i \}, \; i=1,\dots, p, \\
S^-_j&=&\{ -1-\epsilon \leq x\leq -1+\epsilon, \; z=-y+b_j \}, \;
j=1,\dots, q.
\end{eqnarray*}
Note that $\pi(S^+_i)=\pi(k_i)$ and $\pi(S^-_j)=\pi(l_j)$. Let $R_k
\subset P$ be given. Then we can write
\begin{center}
$\partial R_k= C_k^1 \cup C_k^2 \cup C_k^3 \cup C_k^4$ where $C_k^1
\subset \pi(k_{i_1}) $, $C_k^2 \subset \pi(l_j)$, $C_k^3 \subset
\pi(k_{i_2}) $, $C_k^4 \subset \pi(l_{j+1})$
\end{center}
for some $1\leq i_1 < i_2 \leq p$ and $1\leq j \leq q$. Lift
$C_k^1, C_k^2, C_k^3, C_k^4$ (along the $x$-axis) so that the
resulting lifts (which will be denoted by the same letters) are
disjoint Legendrian arcs contained in $k_{i_1}, l_j, k_{i_2},
l_{j+1}$ and sitting on the corresponding strips $S^+_{i_1},
S^-_j, S^+_{i_2}, S^-_{j+1}$. For $l=1, 2, 3, 4$, consider
Legendrian linear arcs $I_k^l$ (parallel to the $x$-axis) running
between the endpoints of $C_k^l$'s as in Figure
\ref{single_rectangle}-(a)\&(b). Along each $I_k^l$ the contact
planes make a $90^{\circ}$ left-twist. Let $B_k^l$ be the narrow
band obtained by following the contact planes along $I_k^l$. Then
define $F_k$ to be the surface constructed by taking the union of
the compact subsets of the above strips (containing corresponding
$C_k^l$'s) with the bands $B_k^l$'s (see Figure
\ref{single_rectangle}-(b)). $C_k^l$'s and $I_k^l$'s together
build a Legendrian unknot $\gamma_k$ in $(\mathbb{R}^3,\xi_0)$,
i.e., we set
$$\gamma_k = C_k^1 \cup I_k^1 \cup C_k^2 \cup I_k^2 \cup C_k^3 \cup I_k^3
\cup C_k^4 \cup I_k^4.$$

\medskip \noindent Note that $\pi(\gamma_k)=\partial R_k$,
$\gamma_k$ sits on the surface $F_k$, and $F_k$ deformation
retracts onto $\gamma_k$. Indeed, by taking all strips and bands
in the construction small enough, we may assume that contact
planes are tangent to the surface $F_k$ only along the core circle
$\gamma_k$. Thus, $F_k$ is the ribbon of $\gamma_k$. Observe that,
topologically, $F_k$ is a positive (left-handed) Hopf band.

\begin{figure}[ht]
  \begin{center}
   \includegraphics{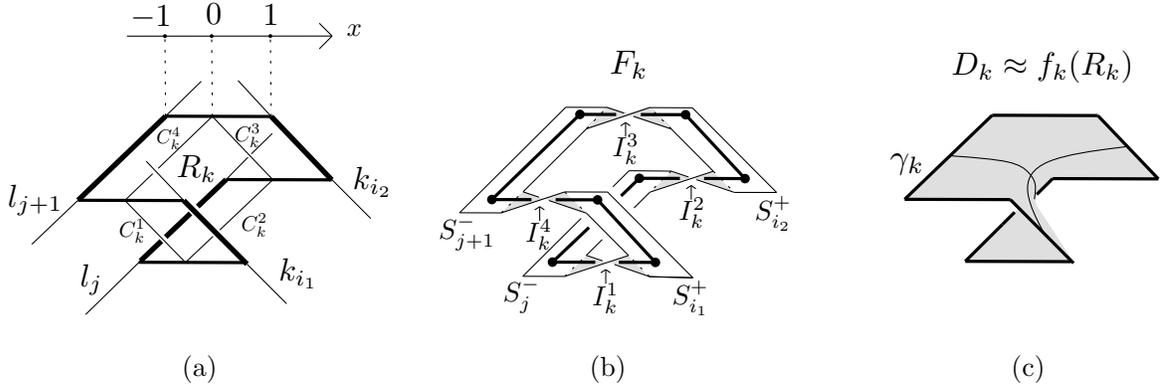}
   \caption{(a) The Legendrian unknot $\gamma_k$, (b) The ribbon $F_k$, (c) The disk $D_k$
  (shaded bands in (b) are the bands $B_k^l$'s) }
  \label{single_rectangle}
    \end{center}
\end{figure}

\noindent Let $f_k : R_k \longrightarrow \mathbb{R}^3$ be a
function modelled by $(a,b) \mapsto c=a^2-b^2$ (for an appropriate
choice of coordinates). The image $f_k(R_k)$ is, topologically, a
disk, and a compact subset of a saddle surface. Deform $f_k(R_k)$
to another ``saddle'' disk $D_k$ such that $\partial D_k =
\gamma_k$ (see Figure \ref{single_rectangle}-(c)). We observe here
that $tw(\gamma_k, D_k)=-1$ because along $\gamma_k$, contact
planes rotate $90^{\circ}$ in the counter-clockwise direction
exactly four times which makes one full left-twist (enough to
count the twists of the ribbon $F_k$ since $F_k$ rotates with the
contact planes along $\gamma_k$ !).

\medskip \noindent We repeat the above process for each rectangle $R_k$ in $P$
and get the set $$\mathfrak{D}=\{ \; D_k \; | \; D_k \approx
f_k(R_k), \; k=1, \dots, N \}$$ consisting of the saddle disks.
Note that by the construction of $\mathfrak{D}$, we have the
property:
\begin{center}
$\hspace{.5cm}$($\ast$) If any two elements of $\mathfrak{D}$
intersect each other, then they must intersect along a
$\hspace{.3cm}$ contractible subset (a contractible union of linear
arcs) of their boundaries.
\end{center}

\noindent For instance, if the corresponding two rectangles (for
two intersecting disks in $\mathfrak{D}$) have only one common
vertex, then those disks intersect each other along the
(contractible) line segment parallel to the $x$-axis which is
projected (by the map $\pi$) onto that vertex.

\medskip \noindent For each $k$, let $D_k'$ be a disk constructed by
perturbing $D_k$ slightly by an isotopy fixing only the boundary
of $D_k$. Therefore, we have
\begin{center}
($\ast \ast$) $\partial D_k = \gamma_k =\partial D_k' \;$, $\quad
int(D_k) \cap int(D_k')=\emptyset \;$, and $\quad tw(\gamma_k,
D_k')=-1=tw(\gamma_k, D_k)$.
\end{center}
\noindent In the following, we will define a sequence $\{ \;
\Delta_k\; | \; k=1, \dots, N \;\}$ of CCD's for $(S^3,\xi_{st})$.
$\Delta_k^1, \Delta_k^2$, and $\Delta_k^3$ will denote the 1-skeleton,
2-skeleton, and 3-skeleton of $\Delta_k$, respectively. First, take
$\Delta_1^1= \gamma_1$, and $\Delta_1^2 =D_1 \cup_{\gamma_1} D_1'$.
By ($\ast \ast$), $\Delta_1$ satisfies the conditions (1) and (2) of
Definition \ref{Cont_Cell_Dec}. By the construction, any pair of
disks $D_k, D_k'$ (together) bounds a Darboux ball (tight 3-cell)
$U_k$ in the tight manifold $(\mathbb{R}^3, \xi_0)$. Therefore, if
we take $\Delta_1^3=U_1 \cup_{\partial} (S^3-U_1)$, we also achieve
the condition (3) in Definition \ref{Cont_Cell_Dec} ( the boundary
union `` $\cup_{\partial}$'' is taken along $\partial U_1 = S^2 =
\partial(S^3-U_1)$ ). Thus, $\Delta_1$ is a CCD for
$(S^3,\xi_{st})$.

\medskip \noindent Inductively, we define $\Delta_k$ from $\Delta_{k-1}$
by setting
\begin{eqnarray*}
\Delta_k^1 &=& \Delta_{k-1}^1 \cup \gamma_k = \gamma_1 \cup \dots \cup \gamma_{k-1} \cup \gamma_k, \\
\Delta_k^2 &=& \Delta_{k-1}^2 \cup D_k \cup_{\gamma_k} D_k' = D_1
\cup_{\gamma_1} D_1'\cup
\dots \cup D_{k-1} \cup_{\gamma_{k-1}} D_{k-1}' \cup D_k \cup_{\gamma_k} D_k',\\
\Delta_k^3 &=& U_1 \cup \dots \cup U_{k-1} \cup U_k \cup_{\partial}
(S^3-U_1 \cup \dots \cup U_{k-1} \cup U_k)
\end{eqnarray*}

\noindent Actually, at each step of the induction, we are applying
Lemma \ref{ExtendingDecomposition} to $\Delta_{k-1}$ to get
$\Delta_k$. We should make several remarks: First, by the
construction of $\gamma_k$'s, the set
$$(\gamma_1 \cup \dots \cup \gamma_{k-1})\cap\gamma_k$$ is a
contractible union of finitely many arcs. Therefore, the union
$\Delta_{k-1}^1 \cup \gamma_k $ should be understood to be a
set-theoretical union (not a topological gluing!) which means that we
are attaching only the (connected) part ($\gamma_k \setminus
\Delta_{k-1}^1$) of $\gamma_k$ to construct the new 1-skeleton
$\Delta_k^1$. In terms of the language of Lemma
\ref{ExtendingDecomposition}, we are setting $I=\Delta_{k-1}^1
\setminus \gamma_k $ and $J=\gamma_k \setminus \Delta_{k-1}^1$.
Secondly, we have to show that $\Delta_k^2 = \Delta_{k-1}^2 \cup D_k
\cup_{\gamma_k} D_k'$ can be realized as the 2-skeleton of a CCD:
Inductively, we can achieve the twisting condition on 2-cells by
using ($\ast \ast$). The fact that any two intersecting 2-cells in
$\Delta_k^2$ intersect each other along some subset of the
1-skeleton $\Delta_k^1$ is guaranteed by the property $(\ast)$ if
they have different index numbers, and guaranteed by $(\ast \ast)$
if they are of the same index. Thirdly, we have to guarantee that
3-cells meet correctly: It is clear that $U_1, \dots, U_k$ meet with
each other along subsets of the 1-skeleton $\Delta_k^1 (\subset
\Delta_k^2)$. Observe that $\partial (U_1 \cup \dots \cup U_k)=S^2$
for any $k=1, \dots, N$ by $(\ast)$ and $(\ast \ast)$. Therefore, we
can always consider the complementary Darboux ball $S^3-U_1 \cup
\dots \cup U_{k-1} \cup U_k$, and glue it to $U_1 \cup \dots \cup
U_k$ along their common boundary 2-sphere. Hence, we have seen that
$\Delta_k$ is a CCD for $(S^3,\xi_{st})$ with Legendrian 1-skeleton
$\Delta_k^1= \gamma_1 \cup \dots \cup \gamma_k$.

\medskip \noindent To understand the ribbon, say $\Sigma_k$, of
$\Delta_k^1$, observe that when we glue the part $\gamma_k \setminus
\Delta_{k-1}^1$ of $\gamma_k$ to $\Delta_{k-1}^1$, actually we are
attaching a 1-handle (whose core interval is $ (\gamma_k \setminus
\Delta_{k-1}^1) \setminus \Sigma_{k-1}$) to the old ribbon
$\Sigma_{k-1}$ (indeed, this corresponds to a positive
stabilization). We choose the 1-handle in such a way that it also
rotates with the contact planes. This is equivalent to extending
$\Sigma_{k-1}$ to a new surface by attaching the missing part (the
part which retracts onto $(\gamma_k \setminus \Delta_{k-1}^1)
\setminus \Sigma_{k-1}$) of $F_k$ given in Figure
\ref{single_rectangle}-(c). The new surface is the ribbon $\Sigma_k$
of the new 1-skeleton $\Delta_k^1$.

\medskip \noindent By taking $k=N$, we get a CCD $\Delta_N$ of
$(S^3,\xi_{st})$. By the construction, $\gamma_k$'s are only
piecewise smooth. We need a smooth embedding of $L$ into the
$1$-skeleton $\Delta_N^1$ (the union of all $\gamma_k$'s). Away
from some small neighborhood of the common corners of $\Delta_N^1$
and $L$ (recall that $L$ had corners before the Legendrian
isotopies), $L$ is smoothly embedded in $\Delta_N^1$. Around any
common corner, we slightly perturb $\Delta_N^1$ using the isotopy
used for smoothing that corner of $L$. This guaranties the smooth
Legendrian embedding of $L$ into the Legendrian graph
$\Delta_N^1=\cup_{k=1}^{N}\gamma_k$. Similarly, any other corner
in $\Delta_N^1$ (which is not in $L$) can be made smooth using an
appropriate Legendrian isotopy.

\medskip \noindent As $L$ is contained in the 1-skeleton $\Delta_N^1$,
$L$ sits (as a smooth Legendrian link) on the ribbon $\Sigma_N$.
Note that during the process we do not change the contact type of
$L$, so the contact (Thurston-Bennequin) framing of $L$ is still
the same as what it was at the beginning. On the other hand,
consider tubular neighborhood $N(L)$ of $L$ in $\Sigma_N$. Being a
subsurface of the ribbon $\Sigma_N$, $N(L)$ is the ribbon of $L$.
By definition, the contact framing of any component of $L$ is the
one coming from the ribbon of that component. Therefore, the
contact framing and the $N(L)$-framing of $L$ are the same. Since
$N(L) \subset \Sigma_N$, the framing which $L$ gets from the
ribbon $\Sigma_N$ is the same as the contact framing of $L$.
Finally, we observe that $\Sigma_N$ is a subsurface of the Seifert
surface $F_{p,q}$ of the torus link (or knot) $T_{p,q}$. To see
this, note that $P$ is contained in the rectangular region, say
$P_{p,q}$, enclosed by the lines $\pi(k_1), \pi(k_p), \pi(l_1),
\pi(l_q)$. Divide $P_{p,q}$ into the rectangular subregions using
the lines $\pi(k_i), \pi(l_j)$, $i=1, \dots, p$, $j=1, \dots, q$.
Note that there are exactly $pq$ rectangles in the division. If we
repeat the above process using this division of $P_{p,q}$, we get
another CCD for $(S^3,\xi_{st})$ with the ribbon $F_{p,q}$.
Clearly, $F_{p,q}$ contains our ribbon $\Sigma_N$ as a subsurface
(indeed, there are extra bands and parts of strips in $F_{p,q}$
which are not in $\Sigma_N$).

\medskip \noindent Thus, (1), (2) and (3) of the theorem are proved
once we set $\Delta=\Delta_N$, (and so $G=\Delta_N^1$,
$F=\Sigma_N$). To prove (4), recall that we are assuming $p>3, \;
q>3$. Then consider
\begin{center}
$\kappa \doteq$ total number of intersection points of all
$\pi(l_j)$'s with all $h_i$'s.
\end{center}
That is, we define $\kappa \doteq | \{\pi(l_j)\; |\; j=1,\dots, q
\} \cap \{h_i \; |\; i=1,\dots, p \} \; |$. Notice that $\kappa$
is the number of bands used in the construction of the ribbon $F$,
and also that if $\D$ (so $P$) is not a single rectangle
(equivalently $p>2$, $q>2$), then $\kappa<pq$. Since there are
$p+q$ disks in $F$, we compute the Euler characteristic and genus
of $F$ as
$$\chi(F)=p+q - \kappa =2-2g(F)- |\partial F| \Longrightarrow g(F)=\displaystyle{\frac{2
-p-q}{2}}+ \displaystyle{\frac{\kappa}{2}}-
\displaystyle{\frac{|\partial F|}{2}}.$$ Similarly, there are $p+q$
disks and $pq$ bands in $F_{p,q}$, so we get
$$\chi(F_{p,q})=p+q - pq =2-2g(F_{p,q})- |\partial F_{p,q}| \Longrightarrow
g(F_{p,q})=\displaystyle{\frac{2-p-q}{2}}+
\displaystyle{\frac{pq}{2}}- \displaystyle{\frac{|\partial
F_{p,q}|}{2}}.$$ Observe that $|\partial F_{p,q}|$ divides the
greatest common divisor $gcd(p,q)$ of $p$ and $q$, so $$|\partial
F_{p,q}| \leq gcd(p,q) \leq p \Longrightarrow g(F_{p,q}) \geq
\displaystyle{\frac{2-p-q}{2}}+ \displaystyle{\frac{pq}{2}}-
\displaystyle{\frac{p}{2}}.$$ Therefore, to conclude
$g(F)<g(F_{p,q})$, it suffices to show that $pq-\kappa > p-|\partial
F|$. To show the latter, we will show $pq-\kappa-p\geq0$ (this will
be enough since $|\partial F|\neq0$).

\medskip \noindent Observe that $pq-\kappa$ is the number of bands
(along $x$-axis) in $F_{p,q}$ which we omit to get the ribbon $F$.
Therefore, we need to see that at least $p$ bands are omitted in
the construction of $F$: The set of all bands (along $x$-axis) in
$F_{p,q}$ corresponds to the set
$$\{\pi(l_j)\; |\; j=1,\dots, q \} \cap \{\pi(k_i)\; |\; i=1,\dots,
p \}.$$ Notice that while constructing $F$ we omit at least $2$
bands corresponding to the intersections of the lines $\pi(k_1),
\pi(k_p)$ with the family $\{\pi(l_j)\; |\; j=1,\dots, q \}$ (in
some cases, one of these bands might correspond to the
intersection of the lines $\pi(k_2)$ or $\pi(k_{p-1})$ with
$\pi(l_1)$ or $\pi(l_q)$, but the following argument still works
because in such a case we can omit at least $2$ bands
corresponding to two points on $\pi(k_2)$ or $\pi(k_{p-1})$). For
the remaining $p-2$ line segments $h_2, \dots, h_{p-1}$, there are
two cases: Either each $h_i$, for $i=2, \dots, p-1$ has at least
one endpoint contained on a line other than $\pi(l_1)$ or
$\pi(l_q)$, or there exists a unique $h_i, 1<i<p$, such that its
endpoints are on $\pi(l_1)$ and $\pi(l_q)$ (such an $h_i$ must be
unique since no two $v_j$'s are collinear !). If the first holds,
then that endpoint corresponds to the intersection of $h_i$ with
$\pi(l_j)$ for some $j\neq1,q$. Then the band corresponding to
either $\pi(k_i) \cap \pi(l_{j-1})$ or $\pi(k_i) \cap
\pi(l_{j+1})$ is omitted in the construction of $F$ (recall how we
divide $P$ into rectangular regions). If the second holds, then
there is at least one line segment $h_{i'}$, which belongs to the
same component of $L$ containing $h_i$, such that we omit at least
$2$ points on $\pi(k_{i'})$ (this is true again since no two
$v_j$'s are collinear). Hence, in any case, we omit at least $p$
bands from $F_{p,q}$ to get $F$. This completes the proof of
Theorem \ref{existence_of_CCD}.
\end{proof}


\begin{corollary} \label{main_corollary}
Given $L$ and $F_{p,q}$ as in Theorem \ref{existence_of_CCD},
there exists an open book decomposition $\mathcal{OB}$ of
$(S^3,\xi_{st})$ such that
\begin{enumerate}
\item $L$ lies (as a Legendrian link) on a page $F$ of $\mathcal{OB}$,
\item The page $F$ is a subsurface of
$F_{p,q}$
\item The page framing of $L$ coming from $F$ is equal to its contact framing
$tb(L)$,
\item If $p>3$ and $q>3$, then $g(F)$ is
strictly less than $g(F_{p,q})$,
\item The monodromy $h$ of $\mathcal{OB}$ is given by $h=t_{\gamma_1} \circ
\dots \circ t_{\gamma_N}$ where $\gamma_k$ is the Legendrian
unknot constructed in the proof of Theorem \ref{existence_of_CCD},
and $t_{\gamma_k}$ denotes the positive (right-handed) Dehn twist
along $\gamma_k$.
\end{enumerate}
\end{corollary}
\begin{proof} The proofs of (1), (2), (3), and (4) immediately follow from Theorem
\ref{existence_of_CCD} and Lemma \ref{RibbonIsPage}. To prove (5),
observe that by adding the missing part of each $\gamma_k$ to the
previous 1-skeleton, and by extending the previous ribbon by
attaching the ribbon of the missing part of $\gamma_k$ (which is
topologically a 1-handle), we actually positively stabilize the
old ribbon with the positive Hopf band $(H^+,t_{\gamma_k})$.
Therefore, (5) follows.
\end{proof}

\begin{figure}[ht]
  \begin{center}
   \includegraphics{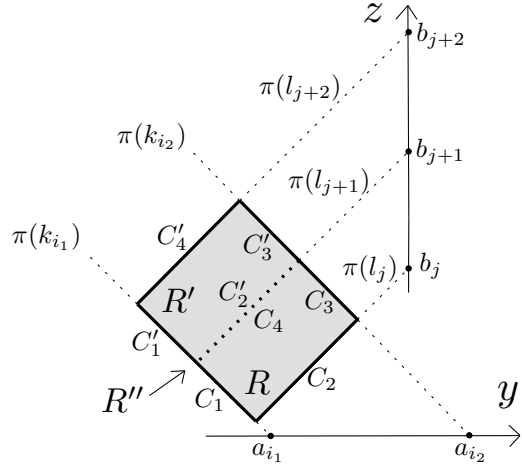}
   \caption{Replacing $R$, $R'$ with their union $R''$}
  \label{combining_rectangles}
    \end{center}
\end{figure}

\noindent With a little more care, sometimes we can decrease the
number of 2-cells in the final 2-skeleton. Also the algorithm can
be modified for split links:

\begin{remark} \label{increasing_efficiency}
Under the notation used in the proof of Theorem
\ref{existence_of_CCD}, we have the following:
\begin{enumerate}

\item Suppose that the link $L$ is split (so $P$ has at least two
connected components). Then we can modify the above algorithm so
that Theorem \ref{existence_of_CCD} still holds.
\item Let $T_j$ denote the row (or set) of rectangles (or
elements) in $P$ (or in $\mathfrak{P}$) with bottom edges lying on
the fixed line $\pi(l_j)$. Consider two consecutive rows $T_j,
T_{j+1}$ lying between the lines $\pi(l_j), \pi(l_{j+1})$, and
$\pi(l_{j+2})$. Let $R \in T_j$ and $R' \subset T_{j+1}$ be two
rectangles in $P$ with boundaries given as
$$\partial R=C_1 \cup C_2\cup C_3 \cup C_4, \quad
\partial R'=C_1' \cup C_2'\cup C_3' \cup C_4'$$
Suppose that $R$ and $R'$ have one common boundary component lying
on $\pi(l_{j+1})$, and two of the other components lie on the same
lines $\pi(k_{i_1}), \pi(k_{i_2})$ as in Figure
\ref{combining_rectangles}. Let $\gamma, \gamma' \subset \Delta_N^1$
and $D, D' \subset \Delta_N$ be the corresponding Legendrian unknots
and 2-cells of the CCD $\Delta_N$ coming from $R, R'$. That is,
\begin{center}
$\partial D=\gamma$, $\partial D'=\gamma'$, and $\pi(D)=R$,
$\pi(D')=R'$
\end{center}
Suppose also that $L \cap \gamma \cap \gamma'=\emptyset$. Then in
the construction of $\Delta_N$, we can replace $R, R' \subset P$
with a single rectangle $R''=R \cup R'$. Equivalently, we can take
out $\gamma \cap \gamma'$ from $\Delta_N^1$, and replace $D, D'$ by
a single saddle disk $D''$ with $\partial D''= (\gamma \cup \gamma')
\setminus (\gamma \cap \gamma')$.
\end{enumerate}
\end{remark}

\begin{proof}
To prove each statement, we need to show that CCD structure and all
the conclusions in Theorem \ref{existence_of_CCD} are preserved
after changing $\Delta_N$ the way described in the statement.

\medskip \noindent To prove (1), let $P^{(1)}, \dots, P^{(m)}$ be the
separate components of $P$. After putting the corresponding separate
components of $L$ into appropriate positions (without changing their
contact type) in $(\mathbb{R}^3,\xi_0)$, we may assume that the
projection
$$P=P^{(1)} \cup \dots \cup P^{(m)}$$
of $L$ onto the second quadrant of the $yz$-plane is given similar as
the one which we illustrated in Figure \ref{L_separate}.
\begin{figure}[ht]
  \begin{center}
   \includegraphics{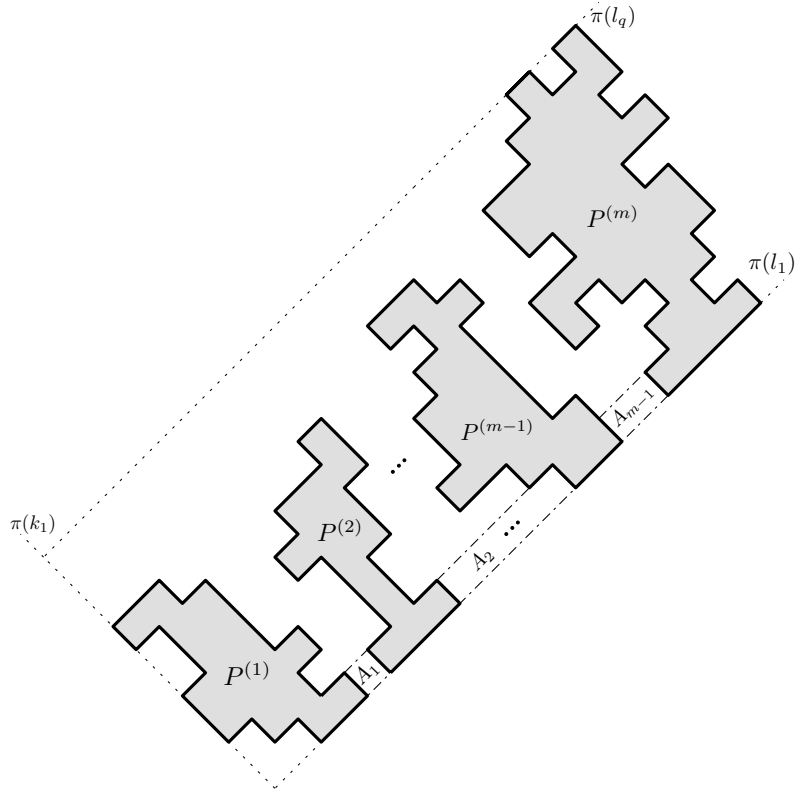}
   \caption{Modifying the algorithm for the case when $L$ is split}
  \label{L_separate}
    \end{center}
\end{figure}

\noindent In such a projection, we require two important properties:
\begin{enumerate}
\item $P^{(1)}, \dots, P^{(m)}$ are located from left to right in
the given order in the region bounded by the lines $\pi(k_1)$,
$\pi(l_1)$, and $\pi(l_q)$.
\item Each of $P^{(1)}, \dots, P^{(m)}$ has at least one edge on the
line $\pi(l_1)$.
\end{enumerate}

\noindent If the components $P^{(1)} \dots P^{(m)}$ remain
separate, then our construction in Theorem \ref{existence_of_CCD}
cannot work (the complement of the union of 3-cells corresponding
to the rectangles in $P$ would not be a Darboux ball; it would be
a genus $m$ handle body). So we have to make sure that any
component $P^{(l)}$ is connected to the some other via some bridge
consisting of rectangles. We choose only one rectangle for each
bridge as follows: Let $A_l$ be the rectangle in $T_1$ (the row
between $\pi(l_1)$ and $\pi(l_2)$) connecting $P^{(l)}$ to
$P^{(l+1)}$ for $l=1, \dots, m-1$ (see Figure \ref{L_separate}).
Now, by adding 2- and 3-cells (corresponding to $A_1, \dots,
A_{m-1}$), we can extend the CCD $\Delta_N$ to get another CCD for
$(S^3,\xi_{st})$. Therefore, we have modified our construction
when $L$ is split.

\medskip \noindent To prove (2), if we replace $D''$ in the way
described above, then by the construction of $\Delta_N^3$, we also
replace two 3-cells with a single 3-cell whose boundary is the union
of $D''$ and its isotopic copy. This alteration of $\Delta_N^3$
does not change the fact that the boundary of the union of all
3-cells coming from all pairs of saddle disks is still homeomorphic
to a 2-sphere $S^2$, Therefore, we can still complete this union to
$S^3$ by gluing a complementary Darboux ball. Thus, we still have a
CCD. Note that $\gamma \cap \gamma'$ is taken away from the 1-skeleton.
However, since $L \cap \gamma \cap \gamma'=\emptyset$, the new
1-skeleton still contains $L$. Observe also that this process does
not change the ribbon $N(L)$ of $L$. Hence, the same conclusions in
Theorem \ref{existence_of_CCD} are satisfied by the new CCD.
\end{proof}


\section{Examples} \label{examples}

\begin{figure}[ht]
  \begin{center}
   \includegraphics{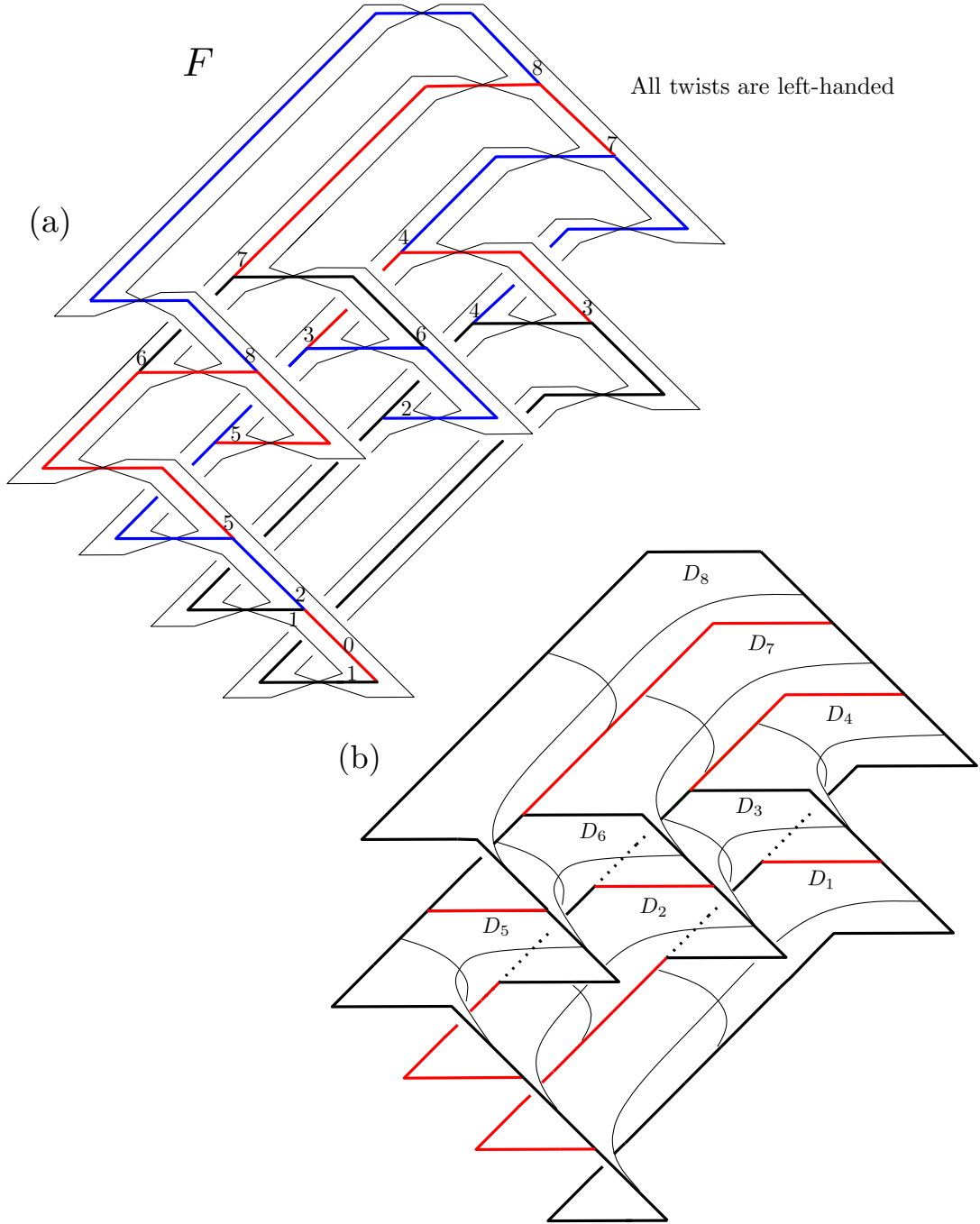}
   \caption{(a) The page $F$ for the right trefoil knot, (b) Construction of $\Delta$}
  \label{trefoil_on_page}
    \end{center}
\end{figure}

\medskip \noindent
\textbf{Example I.} As the first example, let us finish the one
which we have already started in the previous section. Consider
the Legendrian right trefoil knot $L$ (Figure \ref{trefoil_0}) and the
corresponding region $P$ given in Figure \ref{trefoil_2}. Then we
construct the 1-skeleton, the saddle disks, and the ribbon of the
CCD $\Delta$ as in Figure \ref{trefoil_on_page}.

\noindent In Figure \ref{trefoil_on_page}-(a), we show how to
construct the 1-skeleton $G=\Delta^1$ of $\Delta$ starting from a
single Legendrian arc (labelled by the number `` 0 ''). We add
Legendrian arcs labelled by the pairs of numbers ``$1, 1$''$,
\dots, $``$8, 8$'' to the picture one by one (in this order). Each
pair determines the endpoints of the corresponding arc. These arcs
represent the cores of the 1-handles building the page $F$ (the
ribbon of $G$) of the corresponding open book $\mathcal{OB}$. Note
that by attaching each 1-handle, we (positively) stabilize the
previous ribbon by the positive Hopf band $(H^+_k, t_{\gamma_k})$
where $\gamma_k$ is the boundary of the saddle disk $D_k$ as
before. Therefore, the monodromy $h$ of $\mathcal{OB}$ supporting
$(S^3,\xi_{st})$ is given by
$$h=t_{\gamma_1} \circ \dots \circ t_{\gamma_8}$$
where $t_{\gamma_k} \in Aut(F,\partial F)$ denotes the positive
(right-handed) Dehn twist along $\gamma_k$. To compute the genus
$g_F$ of $F$, observe that $F$ is constructed by attaching eight
1-handles (bands) to a disk, and $|\partial F|=3$ where $|\partial
F|$ is the number of boundary components of $F$. Therefore,
$$\chi(F)=1-8=2-2g_F- |\partial F| \Longrightarrow g_F=3.$$

\begin{figure}[ht]
  \begin{center}
   \includegraphics{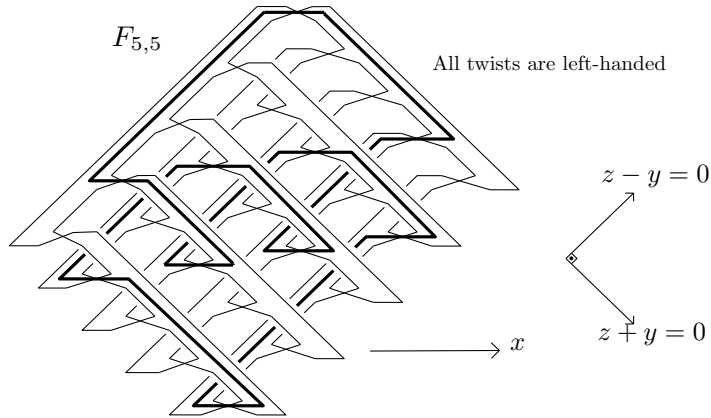}
   \caption{Legendrian right trefoil knot sitting on $F_{5,5}$}
  \label{trefoil_on_5_5_torus}
    \end{center}
\end{figure}

\noindent Now suppose that $(M_1^{\pm},\xi^{\pm}_1)$ is obtained
by performing contact $(\pm 1)$-surgery on $L$. Clearly, the
trefoil knot $L$ sits as a Legendrian curve on $F$ by our
construction, so by Theorem \ref{surgeryonbook}, we get the open
book $(F,h_1)$ supporting $\xi$ with monodromy
$$h_1=t_{\gamma_1} \circ \dots \circ t_{\gamma_8} \circ t_{L}^{\mp 1} \in Aut(F,\partial
F).$$ Hence, we get an upper bound for the support genus invariant
of $\xi_1$, namely, $$sg(\xi_1) \leq 3=g_F.$$ We note that the
upper bound, which we can get for this particular case, from
\cite{AO} and \cite{St} is $6$ where the page of the open book is
the Seifert surface $F_{5,5}$ of the $(5,5)$-torus link (see
Figure \ref{trefoil_on_5_5_torus}).

\begin{figure}[ht]
  \begin{center}
   \includegraphics{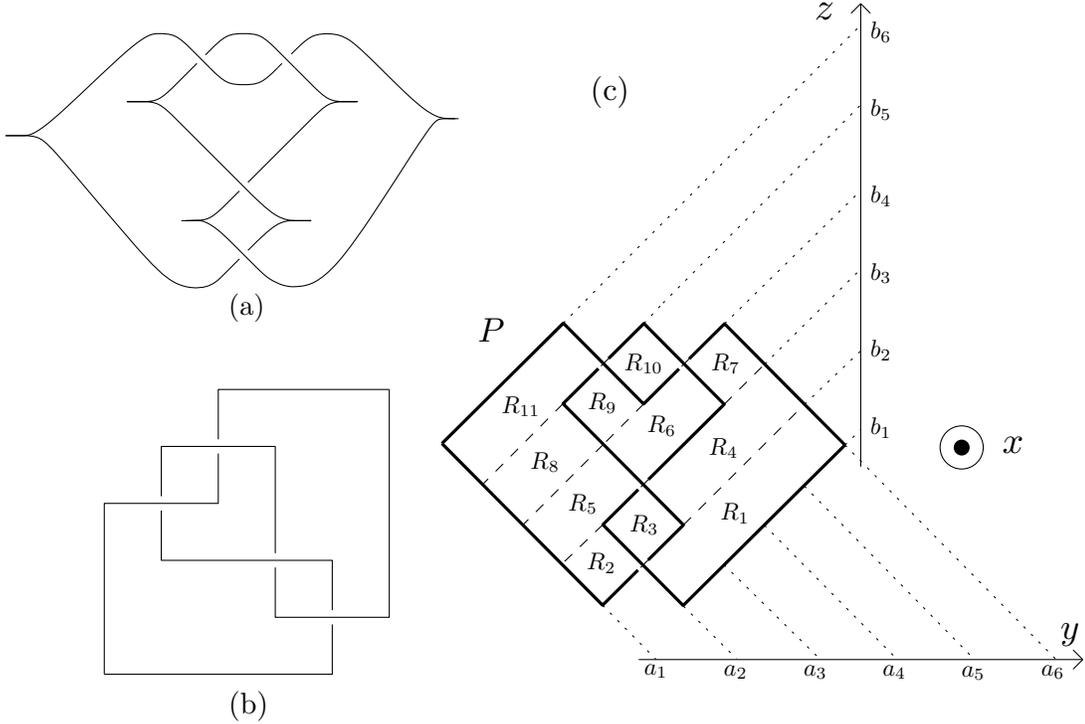}
   \caption{(a),(b) Legendrian figure-eight knot, (c) The region $P$}
  \label{figure8}
    \end{center}
\end{figure}

\begin{figure}[ht]
  \begin{center}
   \includegraphics{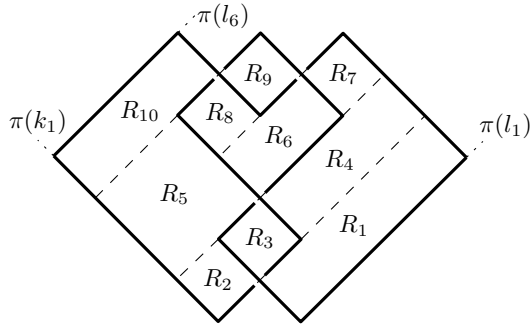}
   \caption{Modifying the region $P$}
  \label{figure8_3}
    \end{center}
\end{figure}

\medskip \noindent
\textbf{Example II.} Consider the Legendrian figure-eight knot
$L$, and its square bridge position given in Figure
\ref{figure8}-(a) and (b). We get the corresponding region $P$ in
Figure \ref{figure8}-(c). Using Remark \ref{increasing_efficiency}
we replace $R_5$ and $R_8$ with a single saddle disk. So this
changes the set $\mathfrak{P}$. Reindexing the rectangles in
$\mathfrak{P}$, we get the decomposition in Figure \ref{figure8_3}
which will be used to construct the CCD $\Delta$.

\begin{figure}[ht]
  \begin{center}
   \includegraphics{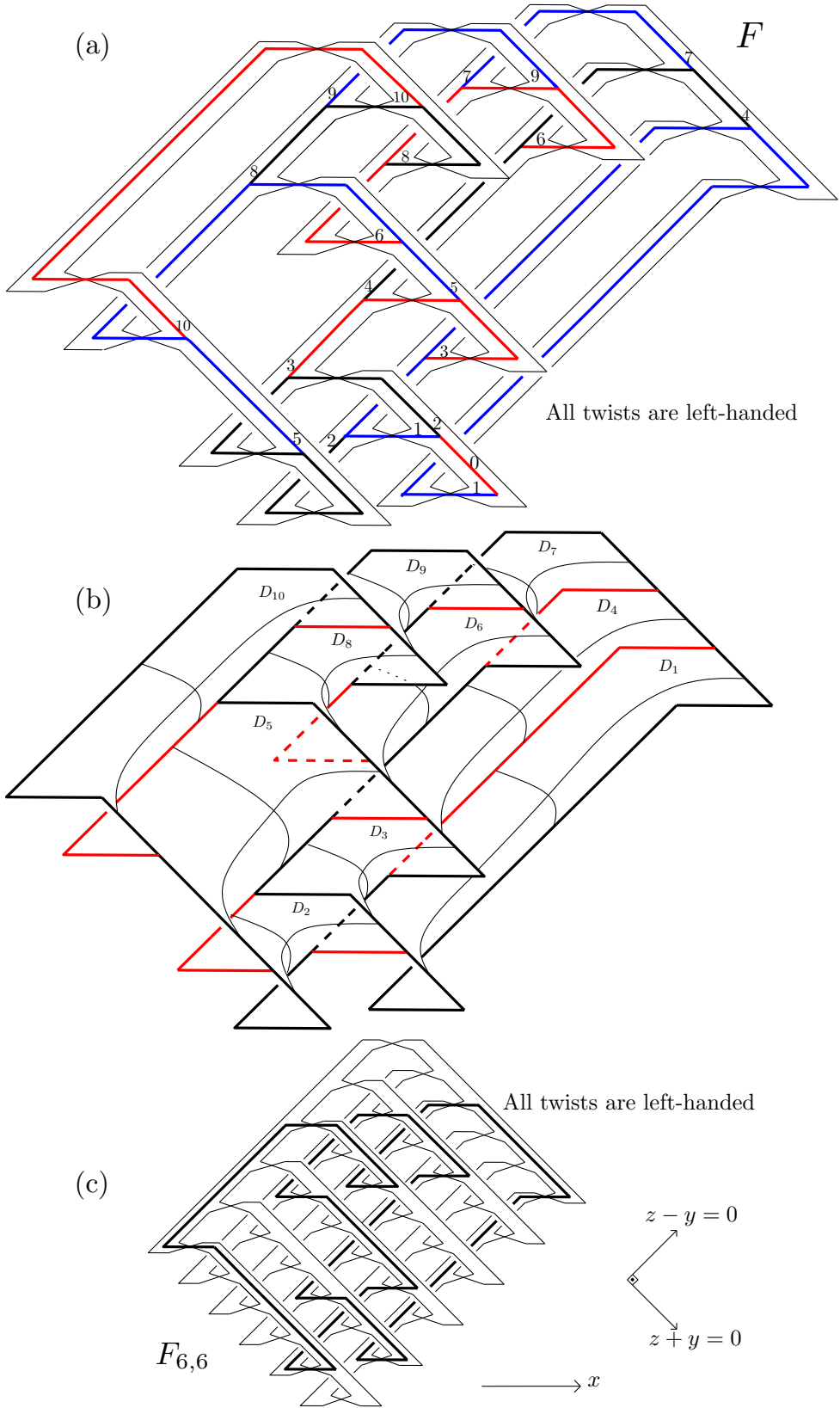}
   \caption{(a) The page $F$, (b) Construction of $\Delta$, (c) The
   figure-eight knot on $F_{6,6}$}
  \label{figure8_on_page}
    \end{center}
\end{figure}

\noindent In Figure \ref{figure8_on_page}-(a), similar to Example
I, we construct the 1-skeleton $G=\Delta^1$ of $\Delta$ again by
attaching Legendrian arcs (labelled by the pairs of numbers ``$1,
1$''$, \dots,$ ``$10, 10$'') to the initial arc (labelled by the
number ``0'') in the given order. Again each pair determines the
endpoints of the corresponding arc, and the cores of the 1-handles
building the page $F$ (of the corresponding open book
$\mathcal{OB}$). Once again attaching each 1-handle is equivalent
to (positively) stabilizing the previous ribbon by the positive
Hopf band $(H^+_k, t_{\gamma_k})$ for $k=1, \dots, 10$. Therefore,
the monodromy $h$ of $\mathcal{OB}$ supporting $(S^3,\xi_{st})$ is
given by
$$h=t_{\gamma_1} \circ \dots \circ t_{\gamma_{10}}$$

\noindent To compute the genus $g_F$ of $F$, observe that $F$ is constructed
by attaching ten 1-handles (bands) to a disk, and $|\partial
F|=5$. Therefore, $$\chi(F)=1-10=2-2g_F- |\partial F| \Longrightarrow g_F=3.$$

\clearpage

\clearpage

\noindent Let $(M_2^{\pm},\xi^{\pm}_2)$ be a contact manifold
obtained by performing contact $(\pm)$-surgery on the figure-8
knot $L$. Since $L$ sits as a Legendrian curve on $F$ by our
construction, Theorem \ref{surgeryonbook} gives an open book
$(F,h_2)$ supporting $\xi_2$ with monodromy
$$h_2=t_{\gamma_1} \circ \dots \circ t_{\gamma_{10}} \circ t_{L}^{\mp 1} \in Aut(F,\partial
F).$$ Therefore, we get the upper bound $ sg(\xi_2) \leq 3=g_F$.
Once again we note that the smallest possible upper bound, which
we can get for this particular case, using the method of \cite{AO}
and \cite{St} is $10$ where the page of the open book is the
Seifert surface $F_{6,6}$ of the $(6,6)$-torus link (see Figure
\ref{figure8_on_page}-(c)).



\begin{thebibliography}{99999}

\bibitem[AO]{AO} S. Akbulut, B. Ozbagci,
{\em Lefschetz fibrations on compact Stein surfaces}, Geom. Topol.
{\bf 5} (2001), 319--334 (electronic).

\bibitem[DG1]{DG1}
F. Ding and H. Geiges, {\em A Legendrian surgery presentation of
contact $3$-manifolds}, Math. Proc. Cambridge Philos. Soc. 136
(2004), no. 3, 583--598.

\bibitem[Et1]{Et1}
J. B. Etnyre, {\em Planar open book decompositions and contact
structures,} IMRN {\bf 79} (2004), 4255--4267.

\bibitem[Et2]{Et2}
J. B. Etnyre, {\em Lectures on open book decompositions and
contact structures}, Floer homology, gauge theory, and
low-dimensional topology, 103--141, Clay Math. Proc., 5, Amer.
Math. Soc., Providence, RI, 2006.

\bibitem[Et3]{Et3}
J.\ Etnyre, \textit{Introductory Lectures on Contact Geometry},
Topology and geometry of manifolds (Athens, GA, 2001), 81--107,
Proc.\ Sympos.\ Pure Math., 71, Amer.\ Math.\ Soc., Providence, RI,
2003.

\bibitem[EO]{EO}
J. Etnyre, and B. Ozbagci, {\em Invariants of Contact Structures
from Open Books}, arXiv:math.GT/0605441, preprint 2006.

\bibitem[Ga]{Ga}
D. Gabai, {\em Detecting fibred links in $S^3$}, Comment. Math.
Helv., 61(4):519-555, 1986.

\bibitem[Ge]{Ge}
H. Geiges, {\em Contact geometry},  Handbook of differential
geometry. Vol. II, 315--382, Elsevier/North-Holland, Amsterdam,
2006.

\bibitem[Gd]{Gd}N Goodman, {\em Contact Structures and Open Books}, PhD thesis, University
of Texas at Austin (2003)

\bibitem[Gi]{Gi} E. Giroux,
{\em G\'eom\'etrie de contact: de la dimension trois vers les
dimensions sup\'erieures}, Proceedings of the ICM, Beijing 2002,
vol. 2, 405--414.

\bibitem[Gm]{Gm} R.\, E. Gompf,
{\em Handlebody construction of Stein surfaces,} Ann. of Math. {\bf
148} (1998), 619--693.

\bibitem[GS]{GS} R.\,E. Gompf, A.\,I. Stipsicz,
{\it 4-manifolds and Kirby calculus}, Graduate Studies in Math. {\bf
20}, Amer. Math. Soc., Providence, RI, 1999.

\bibitem[Ho]{Ho} K. Honda,
{\em On the classification of tight contact structures -I}, Geom.
Topol. {\bf 4} (2000), 309--368 (electronic).

\bibitem[LP]{LP} A. Loi, R. Piergallini,
{\em Compact Stein surfaces with boundary as branched covers of
$B^4$,} Invent. Math. {\bf 143} (2001), 325--348.

\bibitem[Ly]{Ly} H. Lyon,
\textit{Torus knots in the complements of links and surfaces},
Michigan Math. J. 27 (1980), 39-46.

\bibitem[OS]{OS} B. Ozbagci, A.\,I. Stipsicz,
{\it Surgery on contact 3-manifolds and Stein surfaces}, Bolyai
Society Mathematical Studies, {\bf 13} (2004), Springer-Verlag,
Berlin.

\bibitem[Pl]{Pl} O. Plamenevskaya,
\textit{Contact structures with distinct Heegaard Floer invariants},
Math. Res. Lett., \textbf{11} (2004), 547-561.

\bibitem[St]{St} A.\,I. Stipsicz,
{\it Surgery diagrams and open book decomposition of contact
3-manifolds}, Acta Math. Hungar, {\bf 108 (1-2)} (2005), 71-86.

\bibitem[TW]{TW} W. ~P. Thurston, H.~E. Winkelnkemper,
{\em On the existence of contact forms}, Proc. Amer. Math. Soc.
{\bf 52} (1975), 345--347.

\end{thebibliography}
\end{document}